\documentclass{sig-alternate}

\usepackage{color}

\newtheorem{thm}{Theorem}[section]

\newtheorem{definition}[thm]{Definition}

\numberwithin{equation}{section}

\newcommand{\bR}{{\mathbb{R}}}

  \newcommand{\F}{{\mathcal{F}}}

\renewcommand{\P}{{\mathcal{P}}}

\renewcommand{\S}{{\mathcal{S}}}

\begin{document}
%
\conferenceinfo{SAC'09}{March 8-12, 2009, Honolulu, Hawaii,
U.S.A.}

\CopyrightYear{2009} 
      \crdata{978-1-60558-166-8/09/03}  

\title{
Infinite Bar-Joint Frameworks}

\numberofauthors{2}

\author{
\alignauthor \color{black}{J.C. Owen}\\
       \affaddr{D-Cubed, Siemens PLM Software Ltd.,}\\
       \affaddr{Park House, Castle Park,}\\
       \affaddr{Cambridge, U.K.}\\
       \email{owen.john.ext@siemens.com}
\alignauthor \color{black} S.C. Power \\
       \affaddr{Department of Mathematics and Statistics}\\
       \affaddr{Lancaster University}\\
       \affaddr{Lancaster, U.K.}\\
       \email{s.power@lancaster.ac.uk}
}
\date{7 Aug 2008}
\maketitle

\begin{abstract} Some aspects of a mathematical theory
of rigidity and flexibility are developed for general infinite
frameworks and two main results are obtained. In the first
sufficient conditions, of a uniform local nature, are obtained for
the existence of a proper flex of an infinite framework. In the
second it is shown how continuous paths in the plane may be
simulated by infinite Kempe linkages.
\end{abstract}


\category{J.2}{Physical Sciences and Engineering}{Mathematics and
Statistics} \category{J.6}{Computer Aided
Engineering}{Computer-aided design (CAD)}

\terms{Infinite Equation Sets, Materials Analysis, Geometric
Constraints}

\keywords{bar-joint framework, rigidity matrix, rigidity operator,
compactness, Kempe linkage}

\section{Introduction}

We describe some results and work in progress in the analysis of
infinite bar-joint frameworks, their constraint systems and their
solution spaces. In particular we are interested in forms of
flexibility and rigidity.

The behaviour of some physical systems, such as flexible materials
(e.g. foam \cite{des-et-al}) or the positioning of large arrays of
components in an engineering design, may be approximated by a
large number of polynomial equations and effectively modeled by an
infinite equation set. Recently Deshpande et al \cite{des-et-al},
Guest and Hutchinson \cite{gue-hut}, Donev and Torquato
\cite{don-tor} and others have considered rigidity issues for
infinite periodic lattice structures that are of significance for
material analysis. There is a well-established connection between
the rigidity theory of pin-joined structures and geometric
constraint equations  for finite systems
and we expect that a more general study of infinite frameworks
will have relevance for infinite systems of geometric constraints.

 In this article we take a fundamentally
mathematical perspective and lay down a variety of examples,
definitions and concepts to identify some of the diversity of
infinite frameworks per se, whether periodic or not. There are
several motivations for this. The subject in itself, as a
mathematical topic, is novel, intriguing and hybrid, and is able
to draw on established rigidity theory for inspiration and
conjectures. Secondly, there are diverse areas of mathematics that
can be brought to bear or which are appropriate for restricted
classes of frameworks. For example rigidity matrices, being now
infinite, can be viewed as operators on appropriate restricted
spaces of flex vectors, and so are amenable to operator theory
methods and functional analysis perspectives. Thirdly we expect
that the analysis of periodic and aperiodic structures can benefit
by being set in a more general area of analysis which in turn will
lead naturally to the consideration of robust forms of rigidity.

Examples are important. We give several contrasting examples
together with a range of concepts and terminology aimed at
differentiating some of the rich variety of framework types. We
follow this with two main results. The first, with full proof,
illustrates one way in which topological compactness in function
spaces can be useful. Here it assists in establishing sufficient
conditions for the continuous (real) flexibility of an infinite
framework. The second result is inspired by the celebrated 1876
linkage construction of Kempe who showed that a finite linkage (a
two-dimensional bar-joint framework with one degree of
flexibility) can be designed to simulate a given algebraic curve.
See also Gao et al \cite{gao-et-al}. Here we show how infinite
frameworks can simulate  continuous functions, once again, with
zero error, and we provide outline proofs.

For diverse discussions of finite framework rigidity and
constraint systems see, for example, \cite{asi-rot},
\cite{gra-ser-ser}, \cite{lam}, \cite{owe}, \cite{owe-pow-1},
\cite{par-how}, \cite{roth}, \cite{whi-1}.

\section{Examples}

We define a (countable) \emph{infinite} (bar-joint)
\emph{framework} in $\mathbb{R}^{d}$ to be a pair $\mathcal{G}=
(G,p)$ where $G=(V,E)$, the abstract graph of $\mathcal{G}$, has
countable vertex set $V$ and edge set $E$, and where $p=(p_{1},
p_{2}, \dots)$ with $p_{i} \in \mathbb{R}^{d}$ for all $i$, is the
\emph{framework vector} of $\mathcal{G}$ associated with an
enumeration $V=\{v_{1}, v_{2}, \dots\}$.
 The \emph{framework
edges} of $\mathcal{G}$ are the unordered straight line segments
$[p_{i},p_{j}]$ for each $i,j$ with $(v_{i},v_{j})$ an edge of
$G$.

The following example in one dimension is revealing.

Let $(G,p)$ be the infinite framework in $\mathbb{R}$ with
framework vector $p=(p_{1},p_{2}, \dots)$  and  framework edges
$[p_{n},p_{n+1}]$ for all $n$. The abstract graph $G$ here is a
tree with a single branch and a single vertex of degree $1$. Two
such linear frameworks  $(G,p)$ and  $(G,q)$ are equivalent if
$|q_{n}-q_{n+1}|= |p_{n}-p_{n+1}|$, for $n=1,2, \dots$, and are
congruent if for some isometry of $T$ of $\mathbb{R}$, we have
$q_n = Tp_n$ for all $n$. Recall the fact that for every real
number $\alpha \in \mathbb{R}$ there is a sequence $a_{1},a_{2},
\dots $ with $a_{n}=1$ or $-1$ for all $n$, such that $\alpha=
\sum_{n=1}^{\infty} a_{n}n^{-1}$. Thus the framework with vector
$p =(0,1,1-1/2,1-1/2+1/3,...)$, has uncountably many pairwise
noncongruent equivalent frameworks (obtained by flipping edge
directions). From this, and analogous infinitely folding
frameworks in higher dimensions, we also easily see  that a
continuously rigid framework (formally defined below) may possess
uncountably many pairwise noncongruent equivalent frameworks that
are $\epsilon$-close (in the sense of Definition 3.2).

\subsection{Diminishing Rectangles}

Let $\mathcal{G}_{1}= (G,p)$ be the infinite planar framework in
Figure 1. We may label it with $p_0=(1,-1/4)$,  $p_{1}=(1,0),$ $
p_{2}=(1,1), p_{3}=(\frac{1}{2}, 0), p_{4}= (\frac{1}{2},
\frac{1}{2}), p_{5}=(\frac{1}{3}, 0), p_{6}= (\frac{1}{3},
\frac{1}{3})$, and so on, with edges $[p_{i}, p_{i+1}]$, for $i$
odd, with edges $[p_{i}, p_{i+2}]$ for all $i\geq 1$, and with the
indicated edges to $p_0$ which have the effect of "rigidifying"
the $x$-axis edges.

\begin{center}
\begin{figure}[h]
\centering
\includegraphics[width=5cm]{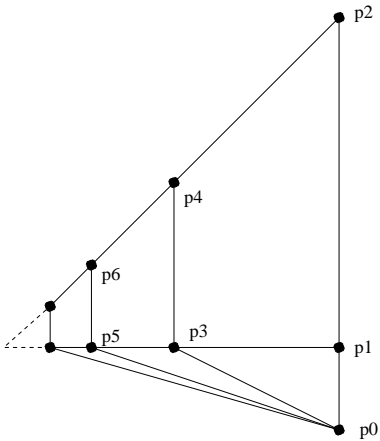}
\caption{An inflexible framework.}
\end{figure}
\end{center}

Suppose for the moment that $p_{i}(t), i=0,1,2, \dots$ are
continuous functions from $[0,1]$ to $\mathbb{R}^{2}$, with
$|p_{i}(t)-p_{j}(t)|= |p_{i}-p_{j}|$ for all $t$ and all framework
edges $[p_{i}, p_{j}]$. We may suppose moreover that $p_{i}(t)=
p_{i}$ for all $t \in [0,1]$ and $i$ odd. Then it can be proven
that $p_{i}(t)= p_{i}$ for all $t$ and all $i$ even. That is, the
\emph{continuous flex} $p(t)= (p_{1}(t), p_{2}(t), \dots )$ must
be constant. The reason for this, roughly speaking, is that the
flexible rectangular subframework determined by $p_{i}, p_{i+1},
p_{i+2}, p_{i+3}$, for $i=1,3,5, \dots$, has a limited
flexibility, tending to zero as $i$  tends to infinity, and since
flexes propogate linearly no continuous flex of $p_{2}$ (and hence
of any $p_{i}$) is admissible.
\bigskip

\subsection{Cobweb Graph Frameworks}

Let $\mathcal{G}_{1}$ be the square frame framework with framework
points
\[ \{p_{1}, \dots, p_{4}\}= \{ (1,1), (1,-1), (-1,-1), (-1,1) \} .\]
Let $\mathcal{G}_{\infty}$ be the framework which, roughly
speaking, consists of the union $\mathcal{G}_{1} \cup \frac{1}{2}
\mathcal{G}_{1} \cup \frac{1}{4} \mathcal{G}_{1} \cup \dots $
together with connected edges between the corresponding corners of
consecutive squares. Then we call $\mathcal{G}_{\infty}$ the
dyadic \emph{cobweb framework} and we have $\mathcal{G}_{\infty}=
(G_{\infty},p)$ where the abstract graph $G_{\infty}$ is a cobweb
graph. It can be shown that while every finite subframework of
$G_{\infty}$ is continuously flexible, $G_{\infty}$ itself is not,
again for reasons of vanishing flexibility, although in this (less
intuitive) case some geometric analysis is needed.

\begin{center}
\begin{figure}[h]
\centering
\includegraphics[width=5cm]{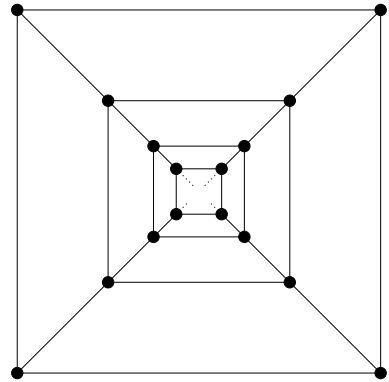}
\caption{The dyadic cobweb framework.}
\end{figure}
\end{center}

The similar framework $G^{\infty}$ which is constructed on the
framework points of \quad $\mathcal{G}_{1} \cup 2\mathcal{G}_{1}
\cup 4\mathcal{G}_{1} \cup \dots$ \quad is continuously flexible,
while the two-way infinite framework
$\mathcal{G}_{\infty}^{\infty}= \mathcal{G}_{\infty} \cup
\mathcal{G}^{\infty}$ is rigid.

From a mathematical perspective (and perhaps also from other
perspectives) the cobweb framework $G^{\infty}$ is interesting in
that it admits a proper flex which is increasingly negligible
towards infinity. We see an opposite amplifying effect in the next
example.

\subsection{Lattice Flexing}

It is straightforward to construct a finite framework with one
degree of flexibility which 'simulates'  two rigid bars jointed at
their midpoints. For example take four equal length framework
edges joined at a common central framework point and add two
"extraneous" vertices and six edges to force them to be colinear
in two  pairs.

Similarly we can simulate two rigid bars jointed at any interior
points. Cocatenating infinitely many such 'tweezer' components
leads to frameworks with one degree of flexibility. Cocatenating
identical components leads to the infinite wine rack in the
diagram.  (The open circles in Figure 3 indicate interior jointing
of rigid bars.)

Note that any proper flex $p(t)= (p_{1}(t),p_{2}(t), \dots )$ of
the infinite winerack  is \emph{unbounded} in the sense that for
each $t
> 0$ the sequence $p(t)-p(0)$ is not a bounded sequence.

\begin{center}
\begin{figure}[h]
\centering
\includegraphics[width=5cm]{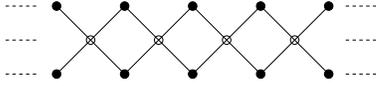}
\caption{The infinite winerack framework.}
\end{figure}
\end{center}

One can assemble infinitely many tweezer components in all manner
of interesting ways. In particular one may arrange the total edge
length sum to be finite while maintining flexibility (despite the
presence of arbitrarily small rectangles). One can also arrange
tree structured assemblages $\mathcal{G}$ with  Cantor set
topological boundaries which exhibit interesting dynamics under
framework flexing. An example of this is the Cantor tree framework
in Figure 4.

\begin{center}
\begin{figure}[h]
\centering
\includegraphics[width=5cm]{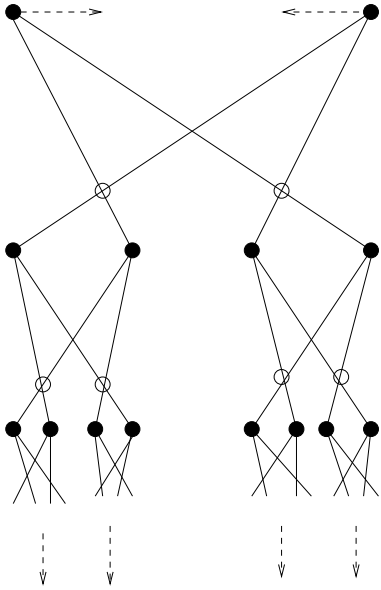}
\caption{Cantor tree tweezer framework.}
\end{figure}
\end{center}

\subsection{Periodic Frameworks}
\vspace{.2in}

 Spatially periodic frameworks are, of course,
ubiquitous, appearing, for example, in the mathematical models
underlying crystallography and polymer frameworks, in the real
finite world of space structures, and in the pure mathematical
realm of planar tilings \cite{gru-she}. Simply enumerating
periodic tetrahedral frameworks (of interest for hypothetical
tetrahedral $SiO_2$) is a major project (for which see Treacy et
al \cite{tre-et-al}). However, as Donev and Torquato and others
have observed there has been little development of rigidity theory
for truly periodic (and hence infinite) frameworks.

We do not comment further on this here except to refer the reader
to Deshpande et al \cite{des-et-al}, Guest and Hutchinson
\cite{gue-hut}, and Donev and Torquato \cite{don-tor}, for
examples of interest in materials analysis, and to remark that
some actual space structures are, in a manner of speaking, almost
infinite. (The dome of the Sports Palace Sant Jordi in Barcelona
was assembled from 9070 bars and 2343 joints \cite{par-how}.)

\section{Rigidity and Rigidity Operators}

\begin{definition} Frameworks $\mathcal{G}=(G,p)$ and
$\mathcal{G}'=(G',p')$ are \emph{equivalent} if there is a graph
isomorphism $\pi : G \rightarrow G'$ such $|p_{i}-p_{j}|=
|p_{\pi(i)}'- p_{\pi(j)}'|$ for all edges $(v_{i}, v_{j})$ of $G$
(where $v_{\pi(i)}'= \pi(v_{i})$). The frameworks are
\emph{congruent} if $ Tp_{i}= p_{\pi(i)}'$ for all $i$ for some
permutation $\pi$ and isometry $T$ of $\mathbb{R}^{d}$.
\end{definition}

For a useful discussions of equivalence in the finite case,
including the problem of unique rigidity (or global rigidity), in
which equivalent frameworks are necessarily congruent, see
Connelly \cite{con}.

\begin{definition} A framework $(G,p)$ is
$\epsilon$-rigid
whenever $(G',p')$ is an equivalent framework (with equivalence
map $\pi = $identity) and $|p_{i}-p_{i}'| \leq \epsilon$ for all
$i$, then $(G,p)$ and $(G',p')$ are congruent. A framework $(G,p)$
is \emph{perturbationally rigid} if it is $\epsilon$-rigid for
some $\epsilon >0$.
\end{definition}

The concept of $\epsilon$-rigidity was introduced in the
pioneering paper of Gluck \cite{glu} for finite frameworks. For
finite frameworks it was shown by Gluck
 to be
equivalent to continuous rigidity, as expressed in the next
definition, and also, in the case of generic frameworks, to
infinitesimal rigidity, as expressed in the subsequent one. (A
generic finite framework is one whose framework point coordinates
are algebraically independent over the rational numbers.)

It is convenient to
 restrict to two-dimensional frameworks.

\begin{definition}
Let  $(G, p)$ be a (possibly infinite) framework in $\bR^2$ with
connected abstract graph  $G=(V, E)$. Let  $V = \{v_1, v_2, \dots
\}$ and $p=(p_1,p_2, \dots )$. Then

(a) $(G, p)$ is said to be \textit{flexible}, or more precisely,
\textit{continuously flexible}, with a \textit{(proper) continuous
flex} $p(t)$ if there exists a function $p(t) = (p_1(t), p_2(t),
\dots )$
 from  $[0,1]$ to $\prod_V \bR^2$ with the following five
 properties.

(i) $p(0) = p$,

 (ii) each coordinate function $p_i(t)$ is  continuous,

(iii) for some base edge $(a,b)$ with $|p_a - p_b| \neq 0$,
$p_a(t) = p_a(0)$ and $p_b(t) = p_b(0)$ for all $t$,

(iv) each edge distance is conserved; $|p_i(t) - p_j(t)| = |p_i(0)
- p_j(0)| $ for all edges $(v_i , v_j)$, and all $t$,


(v) $p(t)$  is not a  constant function.

(b) The framework $(G, p)$ is \textit{rigid} (or continuously
rigid) if it is not flexible, that is, if it has no (proper)
continuous flex.
\end{definition}

We have already seen from our elementary linear examples that
perturbational rigidity may fail rather spectacularly for a
continuously rigid framework. This can also be seen in a similar
way for the simple infinite framework suggested by  Figure 5 (and
this framework is "regular" in the terminology below).

\begin{center}
\begin{figure}[h]
\centering
\includegraphics[width=5cm]{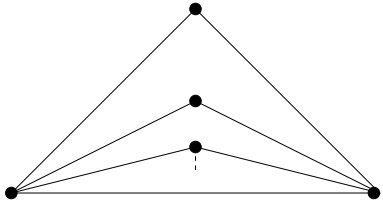}
\caption{Rigid but not $\epsilon$-rigid for all $\epsilon$.}
\end{figure}
\end{center}

If (instantaneous velocity) vectors $u_{1},u_{2}, \dots $ in
$\mathbb{R}^{2}$ have the property $(p_{i}-p_{j}).(u_{i}-u_{j})=
0$ for all $i,j$ then the vector  $u=(u_{1},u_{2}, \dots)$ in the
infinite dimensional vector space $\mathcal{H}_{v}= \mathbb{R}^{2}
\oplus \mathbb{R}^{2} \oplus \dots $ is called an
\emph{infinitesimal flex}. We maintain this traditional
terminology for infinite frameworks even though $u$ may be an
unbounded sequence. Every framework in $\mathbb{R}^{2}$ has a
three-dimensional subspace of infinitesimal flexes coming from
isometric motions (spatial isometries). Any nonzero infinitesimal
flex not in this space is a \textit{proper infinitesimal flex}.

\begin{definition} An infinite framework $(G,p)$ is
\emph{infinitesimally rigid} if it has no proper infinitesimal
flexes and is \emph{infinitesimally flexible} otherwise.
\end{definition}

The dyadic cobweb framework has infinitely many infinitesimal
motions  because of its symmetries (which allow evident "local
infinitesimal rotations"), but even a generic cobweb framework
retains a proper infinitesimal flex. (By generic here we mean
merely that each finite subgraph is generic.) So, in contrast with
the finite case, generic infinite frameworks may be continuously
rigid without being infinitesimally rigid. Also infinitesimal
rigidity and perturbational rigidity differ, so the three
definitions are pairwise inequivalent in general.

We now indicate briefly below how the three definitions also have
conditional forms that are appropriate to infinite frameworks.

The following terminology is useful. Let $\mathcal{G}= (G,p)$ be
an infinite framework with countably many edges and vertices, and
let $e_{1},e_{2}, \dots$ be an enumeration of the edges. Let
$d_{e}=d_{ij}$ be the separation $|p_{i}-p_{j}|$ for the edge
$e=(v_{i},v_{j})$. Say that $\mathcal{G}$ is \emph{edge vanishing}
(respectively \emph{edge unbounded}) if the sequence
$(d_{e_{i}})_{i=1}^{\infty}$ has no lower bound $\delta > 0$
(respectively no upper bound), and say that $\mathcal{G}$ is
\emph{regular} if a lower bound $\delta> 0$ and an upper bound
$M>0$ exist. We also say that $\mathcal{G}$ is \emph{bounded} or
\emph{unbounded} if the sequence $p$ has this property. Also
$\mathcal{G}$ is \emph{locally finite} if each vertex of $V$ has
finite degree. In particular, periodic frameworks, with a
repeating finite cell
(\cite{des-et-al},\cite{don-tor},\cite{gue-hut}), are regular and
locally finite.

In the case of edge vanishing frameworks in many respects it is
appropriate to take into account local scales when considering a
perturbed or nearby framework. To quantify this let $m_{i}=
\inf_{j} d_{ij}, M_{i}= \sup_{j}d_{ij}$.  We say that a locally
finite framework $\mathcal{G}=(G,p)$ is \emph{relatively
$\epsilon$-rigid} if every equivalent framework $\mathcal{G}=
(G,p')$ (with $\pi=$id) such $ |p_{i}'-p_{i}| \leq \epsilon m_{i}$
for all $i$ is congruent to $\mathcal{G}$. It is natural then
(particularly in the light of the simple one-dimensional example
above) to determine conditions for relatively $\epsilon$-rigidity.

Similarly, one can consider conditional forms of infinitesimal
rigidity (resp. continuous rigidity) by restricting attention to
specific subspaces of infinitesimal (resp. continuous) flexes
which, for example, may decay at an appropriate local rate. Or, if
one is concerned with decaying  flexes in a regular framework, one
may impose square summable decay. These and similar perspectives
amount to considering the rigidity matrix (and related matrices of
the framework equation system) as a bounded linear operator,
\textit{the rigidity operator}  between appropriate sequence
spaces.

Recall that the \textit{rigidity matrix} $R(G,p)$ of a framework
$(G,p)$ is
 the Jacobian of the  system of framework edge-length
equations multiplied by $1/2$ and evaluated at the framework
points. We take the same definition for countable frameworks. Thus
rows are indexed by edges and columns by the coordinates $x_i,
y_i$ of $p_i$, $i = 1,2,\dots$. (The entry for the edge
$(v_i,v_j)$ and the coordinate $x_i$ is $(x_i-x_j)$ etc.)


One may consider "conditioning" infinitesimal flexes
$u=(u_1,u_2,\dots )$ by requiring that they lie in the vector
space $\ell^\infty_v$ of bounded sequences. (This rules out
infinitesimal rotations of unbounded frameworks for example.) For
a regular framework the rigidity matrix actually determines a
bounded linear transformation $R(G,p)_{\infty, \infty}$ from
$\ell^\infty_v$ to $\ell^\infty_e$. Moreover we may require
bounded \textit{displacements} of framework points which amounts
to a further restriction on the domain of $R(G,p)$. This applies
in particular to the infinite winerack framework;
the natural "squeeze" infinitesimal flex, while being a bounded
sequence does not give bounded displacements of framework points.
In this sense the winerack is a \textit{boundedly isostatic
framework}.

Finally, note that the framework of Figure 5 while being
infinitesimally rigid and continuously rigid has the "flavour" of
an infinitesimally flexible structure. The following definition
gives a natural notion of approximate flexibility to capture this
and which we expect to lead to a useful form of \emph{strong
rigidity}.

\begin{definition}
A framework $(G,p)$ is approximately flexible if for every
$\epsilon >0$ there is a proper flex $u$ with
\[
|(u_i-u_j).(p_i-p_j)| \leq \epsilon(|u_i|+|u_j|)(|p_i-p_j|)
\]
for all edges $v_i,v_j)$. A framework $(G,p)$ is strongly rigid if
it is not approximately flexible.
\end{definition}


\section{Compactness and Proper Flexes}

If an infinite framework has flexible finite subframeworks then
under what conditions might one conclude the existence of a proper
(continuous) flex ? The dyadic cobweb framework
$\mathcal{G}_\infty$ which is inflexible, with all its finite
subframeworks flexible, shows that some care is needed here. In
this section we give a sample theorem which resolves this
question. It is stated and discussed for planar frameworks but
holds for higher dimensions with the same proof. The proof makes
use of the Ascoli-Arzela compactness theorem in the following
form. A bounded equicontinuous sequence of functions $f_k : [0,1]
\to \bR^n, k=1,2,\dots $ has a convergent subsequence. (See
\cite{rud} or \cite{dav-don} for example.)

\begin{definition}
A continuous flex $p(t)$ of a normalised framework is a
\textit{smooth flex} if each coordinate $p_i(t)$ is differentiable
on $[0,1]$ with continuous derivative $p_i'(t)$, where $p_i'(0)$
and $p_i'(1)$ are right and left derivatives respectively.
Furthermore a smooth flex is a \textit{boundedly smooth flex}, or
\textit{$M$-smooth flex}, if for some $M>0$ and for every pair
$p_i, p_j$ the distance function
\[
d_{ij}(t) = |p_i(t) - p_j(t)|
\]
has bounded derivative, with $|d_{i,j}'(t)| \leq M$ for all $t$ in
$[0,1]$.
\end{definition}

Let $(G,p)$ be an infinite locally finite framework in $\bR^2$
with connected graph and with normalised framework vector $p$, in
the sense that $p_1 = (0,0), p_2=(d_{e_1},0)$. Let us say that a
\textit{standard chain} for $(G,p)$ is any sequence of vertex
induced connected subframeworks $(G_1,p)\subseteq (G_2, p)
\subseteq \dots$ whose union is $(G,p)$. Denote the separation
distance $|p_i-p_j|$ by $d_{i,j}$.

The following theorem, in paraphrase, says that there will be a
proper continuous flex of the infinite framework if there are two
framework points such that every finite framework containing them
has at least one smooth flex which changes the separation of these
points, and these separation changes are bounded away from zero.
In general these smooth flexes need not be related and indeed the
entire framework could have, loosely speaking, many (and even
infinitely many) degrees of freedom.

\begin{thm} Let $(G,p)$ be an infinite locally finite
planar framework with connected
graph, let
$$(G_1,p)\subseteq (G_2, p)
\subseteq \dots ,$$ be a standard chain and let $v_i, v_j$ be
vertices of $G_1$. Suppose that there exist $M>0$ and  $c>0$ and a
sequence of $M$-smooth (normalised)  flexes $p^{(r)}(t)$ of
$(G_r,p)$, for $r=1,2, \dots $, such that for all $r$
\[
|d_{i,j}^{(r)}(1) - d_{i,j}^{(r)}(0)| \geq c.
\]
Then $(G,p)$ is continuously flexible.
\end{thm}

The proof is constructed as an iterated application of the
Ascoli-Arzela theorem and a standard diagonal selection to create
a sequence of  coordinate  functions
\[(q^{k(1,1)}_1(t),q^{k(2,2)}_2(t), \dots )
\]
which (although not infinite flexes) converge (uniformly in
coordinates) to a proper flex $q^*(t)$ as $k(n,n)$ tends to
infinity (with $n$). The inequality ensures that the resulting
limit flex is not trivial.

\begin{proof}
Let $\F_l$ be the set of all $M$-smooth flexes $q:[0,1]\to
\bR^{2|V_l|}$ for $(G_l,p)$. This is a nonempty family of
continuous vector-valued functions which are, moreover,
equicontinuous. Let $q^{(r)}(t), r=1,2, \dots $, be the given
sequence of $M$-smooth flexes. Each of these flexes restricts to a
flex of the first subframework $(G_1,p)$. We can write these
restrictions as $P_1q^{(r)}(t)$ where $P_1$ is the natural
projection from the space of infinite framework vectors to the
space determined by the coordinates for $G_1$. This set of
restrictions is a \textit{bounded} set of equicontinuous
vector-valued functions in $\F_1$. This follows from the
hypotheses on derivatives. By the Arzela-Ascoli theorem there is a
uniformly convergent subsequence, determined by some subsequence
$k(1,n), n=1,2,\dots $ of $k=1,2,\dots$. That is we have obtained
a subsequence $q^{(k(1,n))}(t), n=1,2, \dots $, with the $G_1$
coordinates actually converging to a flex of the subframework
$(G_1,p)$.

 Likewise considering the restrictions
$P_{2}q^k(t)$, for $k=k(1,n), $ $n= 1,2,\dots ,$ there is a
subsequence of this subsequence, say $k(2,n), n=1,2,\dots $ such
that the restrictions $P_{2}q^{k(2,n)}(t)$ converge uniformly to a
continuous flex of $(G_{2},p)$ as $n \to \infty$. Continue in this
manner for the entire standard chain, and select the diagonal
subsequence $k(n,n)$. This has  the property that for each
coordinate location, $m$ say, the sequence of coordinate function
$q^{(k(n,n))}_m(t)$ for $n=1,2,\dots $ converges uniformly to a
continuous function $q^*_m(t)$ as $n\to \infty$. Moreover the
function $q^*(t) = (q^*_1(t), q^*_2(t), \dots )$ is the desired
flex. Note in particular that this limit is a proper flex since
the inequality persists in the limit, that is,
\[
||q^*_i(1)- q^*_j(1)| - |q^*_i(0)-q^*_j(0)|| \geq c.
\]
\hfill
\end{proof}

In fact stronger  forms of this theorem hold. For example it is
enough to require that for $r=1,2,\dots ,$ there are smooth flexes
$p^{(r)}(t)$ of $(G_r,p)$ such that for each $l$ the set of
restrictions of $p^{(r)}(t)$ to $(G_l, p)$, for $r \geq l$, are
uniformly boundedly smooth. This scheme is appropriate for
flexible frameworks similar to or containing an infinite winerack.

\section{Infinite Kempe Linkages}

We state a theorem  due to Kempe \cite{kem} and follow this with a
discussion of exactly what the theorem means and the principal
ideas behind the proof.

\begin{thm}
Every finite algebraic curve in the plane has a linkage
realisation.
\end{thm}

Although Kempe does not define a linkage as a mathematical
construct one may view it, in the spirit of Asimow and Roth
\cite{asi-rot}, as a bar-joint framework whose (normalised)
positions give a real variety which is one dimensional (at regular
points). We take the following more convenient linkesque view
which also serves for infinite frameworks. We let $\langle ~ , ~
\rangle$ denote the usual inner product of real vectors.

\begin{definition}
A plane linkage (resp. infinite plane linkage) is a finite (resp.
infinite) connected framework $\mathcal{G} = (G,p)$ in $\bR^2$
with a degree two vertex $v_1$ with edges $(v_1,v_2)$, $(v_1,v_3)$
and a continuous flex $p(t)$ such that

(i) the cosine angle function
\[
g(t)=\langle p_2(t)-p_1(t), p_3(t)-p_1(t)\rangle
\]
is strictly increasing and

(ii) $p(t)$ is the unique flex $q(t)$ of $\mathcal{G}$ with
$q_i(t)=p_i(t), i=1,2,3$.
\end{definition}

Make a partial normalisation by requiring that $p_1$ and all
$p_1(t)$ are equal to the origin $(0,0)$. We may think of a finite
linkage articulating a motion as the points $p_2(t) ,p_3(t) $ make
changing angles $\theta, \phi,$ respectively, with the $x$ axis.
The framework points move smoothly  if $(\theta, \phi)$ move
smoothly. Identifying a specific "end-point" $p_n$ of the
framework, if  $(\theta, \phi)$ moves smoothly in a
one-dimensional (real) algebraic variety then the endpoint
$p_n(\theta,\phi)$ describes an algebraic curve. In particular,
with $\phi$ fixed, a "circular input" via $\theta$ gives an
algebraic curve $p_n(\theta, \phi)$ with $\theta$ ranging in some
interval. Kempe solved the inverse problem by showing that any
particular finite algebraic curve may be realised as such a
linkage curve for some linkage.

The convenience of the double angle parametrisation comes from the
use of parallellogram and quadrilateral linkages in the assembly
of composite linkages.

Kempe's original  construction (which simulates an algebraic
output curve from a linear input) may be conceived of as a
combination of the following four stages.

1. A parallelogram linkage  $L_1= (R,q)$ with $q_1 $ rooted at the
origin, $q_4$ on a given algebraic curve, provides a $(\theta,
\phi)$ (virtual) curve; $\Phi(\theta, \phi)=0$.

2. The observation that $\Phi(\theta, \phi)=0$ translates into a
multiple angle equation of the form
\[
C= \sum A_n\cos(r_n\theta + s_n\phi + t_n).
\]
Write $f(\theta, \phi)$ for the function given by this finite sum.

3. The construction of a linkage $L_2$ so that for input angles
$\theta, \phi$ the $x$ coordinate of the endpoint $p_n(\theta,
\phi)$ is
 $f(\theta, \phi)$. It is in this stage that Kempe uses an
 assembly argument, indicating basic component
 linkages (translator, multiplier, etc) and how they may be
 combined. See also Gao et al \cite{gao-et-al}.

4.  $L_1, L_2$ are joined together at the origin and their
respective edges, incident to the origin, joined appropriately.
Thus the output angles $(\theta, \phi)$ from $L_1$ become input
angles for $L_2$. As $q_4$ moves on the curve $p_n$ move on the
vertical line $x=C$, and vice versa. (One must also add framework
structure to this join to fix the origin to a "base edge" parallel
to the line $x=C$ and so create a free standing linkage.)

If an infinite linkage $(G,p)$ is such that a subsequence
$p_{n_k}, k=1,2, \dots$ is convergent to $p_*$, say, then the flex
gives rise to a continuous plane curve $p^*(t)$. A sample inverse
result, in the spirit of Kempe's theorem,  is given in the
following. We say that the infinite framework $(G,p)$ is
\textit{pointed} if $p$ is a convergent sequence and if the
sequence of edge lengths tends to zero.

\begin{thm}
Let $f(t), t\in [0,2\pi]$ be a continuous real-valued function
with absolutely summable Fourier series. \\
Then the graph of $f$
has an infinite linkage realisation by a pointed locally finite
linkage $(G,p)$.
\end{thm}

The proof follows a similar format to the breakdown above,
although now the sum is infinite, and some modified assembly
components are needed to ensure that edge lengths diminish to
zero.

A consequence of the theorem is that the motion of limit points of
normalised infinite linkages may fail to be continuously
differentiable in every finite interval.

In fact, more generally, we have found an assembly scheme, based
on uniform approximation rather than Fourier series, which creates
an infinite linkage which realises (with no error) a given
continuous curve. Moreover, if one admits non locally finite
linkages, possessing a single framework point with infinite
degree, then we can arrange that this point coincides with the
curve tracing limit point above. In this way we can obtain the
following theorem. Recall that a continuous planar curve (with
parametrisation) is a continuous function from $[0,1]$ to $\bR^2$.
In particular such a curve can be space filling and so these
mathematical linkages are distinctly curious: with a single input
flex a distinguished framework point may visit \textit{every}
point in a region of positive area !

\begin{thm}
Every continuous planar curve has an infinite linkage realisation.
\end{thm}

\vspace{.2in}



\end{document}